\documentclass[twocolumn]{autart}    

\usepackage{amssymb}
\usepackage[dvips]{epsfig}
\usepackage{amsmath,amssymb}
\usepackage{latexsym}
\usepackage{subfigure}
\usepackage{color}
\usepackage[center]{caption}
\usepackage{graphicx}
\usepackage{url}

\newcommand{\be}{\begin{equation}}
\newcommand{\eq}{\end{equation}}
\newcommand{\ba}{\begin{array}}
\newcommand{\ea}{\end{array}}
\newcommand{\bean}{\begin{eqnarray*}}
\newcommand{\eean}{\end{eqnarray*}}
\newcommand{\bea}{\begin{eqnarray}}
\newcommand{\eea}{\end{eqnarray}}

\newcommand{\R}{\rm {I\kern-2pt R}}

\newcommand{\beq}{\begin{equation}}
\newcommand{\eeq}{\end{equation}}

\newtheorem{theorem}{\bf Theorem}[section]
\newtheorem{remark}{\bf Remark}[section]
\newtheorem{lemma}{\bf Lemma}[section]

\newtheorem{proposition}{\bf Proposition}[section]
\newtheorem{definition}{\bf Definition}[section]

\begin{document}

\begin{frontmatter}

\title{Exponential synchronization of the high-dimensional Kuramoto model with identical oscillators under digraphs} 
\thanks[footnoteinfo]{This work is supported in part by National Natural Science Foundation
(NNSF) of China under Grants 61673012 and a project funded by the Priority
Academic Program Development of Jiangsu Higher Education Institutions (PAPD).
}

\author[]{Jinxing Zhang}\ead{jxingzhang@hotmail.com},  
\author[]{Jiandong Zhu}\ead{zhujiandong@njnu.edu.cn}    

\address{Institute of Mathematics, School of Mathematical Sciences, Nanjing Normal University, Nanjing,
        210023, PRC}  
%
%
%
%
%


\begin{abstract}
For the Kuramoto model and its variations, it is difficult to
analyze the exponential synchronization under the general digraphs due to
the lack of symmetry. 
In this paper, for the high-dimensional Kuramoto model of
identical oscillators, a matrix Riccati differential equation
(MRDE) is proposed to describe the error dynamics. Based on the
MRDE, the exponential synchronization is proved by constructing a
total error function for the case of digraphs admitting spanning
trees. Finally, some numerical simulations are given to illustrate
the obtained theoretical results.
\end{abstract}
\end{frontmatter}

\section{Introduction}
Kuramoto model is one of the most representative mathematical
models of complex dynamical networks, which was first proposed by
Yoshiki Kuramoto to describe and explain the synchronization
phenomena in the real world \cite{kuramoto1975}. Kuramoto model
and its many variations have been applied to many fields such as
neuro-science \cite{Cumin}, power systems \cite{Dorfler}, chemical
engineering \cite{kuramoto1984}, geophysics \cite{Vasudevan2015},
and semiconductor lasers arrays \cite{Kozyreff-PRL-2000}. The
interconnecting network of the original Kuramoto model is a
complete graph or the all-to-all topology. For general
interconnecting topologies, the Kuramoto model composed of $m$
oscillators is described as
 \begin{equation}
\label{eq1}
\dot \theta_{i}=\omega_{i}+k\sum\limits_{j=1}^ma_{ij}\sin(\theta_{j}-\theta_{i}),  \ \ i=1,2,\cdots,m,
\end{equation}
where  $\theta_i$ is the phase of the $i$th oscillator, $\omega_i$
is the natural frequency, $A=(a_{ij})$ is the nonnegative
adjacency matrix of the interconnection graph $\mathcal{G}$ and
$k>0$ is the control gain.
\par
In physics community, the researches pay more attentions on the
thermodynamic limit described by a partial differential equation
as the number the oscillators tends to infinity. For details, we
refer the reader to the surveys \cite{Acebron,Dorfler2}. However,
in system and control community, the main interest lies in the
Kuramoto model composed of finite number of oscillators.
Synchronization is a key issue on Kuromoto model and its
variations. It is said that the {\it phase synchronization} is
achieved if
$$\lim_{t\rightarrow \infty}( \theta_i(t)-\theta_j(t))=0\ \ \forall  \ i,\  j=1,2,\cdots,m.$$
In \cite{Jadbabaie-ACC-2004}, it is shown that, when $\mathcal{G}$
is an undirected graph, the Kuramoto model (\ref{eq1}) can be
rewritten as the compact form as follows:
\begin{equation}
\label{eq_2} \dot \theta=\omega-kB\sin(B^{\mathrm{T}}\theta),
\end{equation}
where $\theta=(\theta_1, \theta_2,\cdots,\theta_m)^{\mathrm{T}}$,
$\omega=(\omega_1, \omega_2,\cdots,\omega_m)^{\mathrm{T}}$, and
$B$ is the incidence matrix of an oriented graph
$\mathcal{G}^{\sigma}$ of $\mathcal{G}$. Just based on the compact
form (\ref{eq_2}), the exponential phase synchronization is proved
for the identical Kuramoto model with undirected graphs
\cite{Jadbabaie-ACC-2004}.  For the case of digraphs,   the
theoretical analysis of  synchronization is relatively difficult.
In \cite{Dong} and \cite{Dorfler}, the contraction property is
adopted to achieve the exponential synchronization for the
Kuramoto model with digraphs.

  An interesting issue  is whether  the synchronization or the exponential
synchronization can be achieved for the high-dimensional Kuramoto
model (Lohe model called in \cite{choi2014}) described as follows:
\begin{equation}
 \label{eq2}
  \dot{r}_{i}=\Omega_{i}r_{i}+k\sum\limits_{j=1}^{m}a_{ij}(r_{j}-\frac{r_{i}^{T}r_{j}}{r_{i}^{T}r_{i}}r_{i}),\ i=1,2,\cdots,m,
\end{equation}
where $r_{i}\in{\mathbf{R}^{n}}$ is the state of the $i$th oscillator,
$\Omega_{i}$ is an $n\times n$ skew-symmetric matrix,
$k>0$ is the control gain, and $A=(a_{ij})\in{\mathbf{R}^{m\times{m}}}$
is the weighted adjacency matrix of the interconnecting network. In \cite{zhu2013},
it has been shown that system (\ref{eq2}) can be reduced to the original Kuramoto
model (\ref{eq1}) as $n=2$,
$$\Omega_i=\left[\begin{array}{cc}
0& -\omega_i\\  \omega_i & 0
\end{array}\right],  r_i=\left[\begin{array}{c}
\cos\theta_i\\ \sin\theta_i
\end{array}\right]. $$
The model (\ref{eq2}) with $\Omega_i=0$ as well as the all-to-all
interconnection is first proposed in \cite{Olfati-Saber-ACC-2006}
as a swarm model on spheres, and is used to solve the max-cut
problem in combinatorial optimization. In
\cite{lohe2009,lohe2010}, some collective dynamical behaviors of
(\ref{eq2}) are shown, which have some potential applications in
quantum synchronization of some quantum devises. For the case of
complete graphs, exponential synchronization is proved by using
the concept of order parameter in \cite{choi2014,chi2014}. In
\cite{zhu2013},\cite{zhu2014}, \cite{Markdahl-IEEECDC-2016} and
\cite{Markdahl-IEEEAC-2018}, the phase synchronization on the unit
hemisphere and the almost global synchronization are investigated
for undirected graphs. For the case of digraphs,
\cite{Lageman-CDC-2016} and \cite{zhang2018} achieve phase
synchronization under some limitations on the initial states. In
\cite{Markdahl-Automatica-2018}, a lifting method is proposed to
analyze the almost globle synchronization on the unit sphere for
digraphs. New exciting development on high-dimensional Kuramoto
model can be seen in \cite{Markdahl2018a,Markdahl2018b}.

However, the theoretical analysis of the exponential
synchronization is much more difficult for general graphs. By
Theorem 13 of \cite{Markdahl-IEEEAC-2018}, the local exponential
synchronization of (\ref{eq2}) with $\Omega_i=0$ is achieved for
connected undirected graphs by using the linearization method.
From Theorem 1 of \cite{Lageman-CDC-2016}, one see that the
exponential synchronization is implemented for weakly connected
and balanced digraphs by using the invariant manifold techniques.
Since a weakly connected and balanced digraph is strongly
connected, the topology  condition for the exponential
synchronization imposed in \cite{Lageman-CDC-2016} is very strong.
To the best of the authors' knowledge, for general digraphs, the
exponential synchronization problem of the high-dimensional
Kuramoto model  is still open.

In this paper, the exponential synchronization  is proved for the
high-dimensional Kuramoto model under a general digraph containing
a spanning tree. A matrix Riccati differential equation is
proposed to describe the dynamics of the synchronization errors
for the first time, which plays an important role in the proof the
exponential synchronization. Finally, numerical simulations are
given to validate the obtained theoretical results.  \par
 The rest of this paper is organized as follows. Section 2 gives some preliminaries and the
 problem statement.
 Section 3 includes our main results. Section 4 shows some simulations.
 Finally, Section 5 is devoted to a summary.\par

\section{Preliminaries and Problem Statement}

 Denote by $\mathcal{G}=(\mathcal{V},\mathcal{E},A)$ the weighted
 digraph of the high-dimensional Kuramoto model (\ref{eq2}), which is
 composed a set of nodes $\mathcal{V}=\{1,2,\cdots,m\}$,
 set of edge $\mathcal{E}\subset{\mathcal{V}\times{\mathcal{V}}}$ and
 a weighted adjacent matrix $A=(a_{ij})$ satisfying
 \begin{equation}
 \label{eq201}
 a_{ij}\begin{cases}
 >0 , \quad (j,i)\in{\mathcal{E}}\\
 =0, \quad {\rm otherwise}.
\end{cases}
\end{equation}
Here, a directed edge $(j,i)$ means that  the state information of
agent $j$ can be transmitted to $i$. A sequence of directed edges
$(k_1,k_2),(k_2,k_3),\cdots,(k_{s-1},k_s)$ is called a {\it path}
from $k_1$ to $k_s$, denoted by $k_1k_2\cdots k_s$.
If ${\mathcal G}$ has a node $v_0$ such that, for any another node
$i$, there exists a directed path from $v_0$ to $i$, then $v_0$ is
called a {\it root node} of ${\mathcal G}$. It is well-known that
a digraph has a root node if and only if it has a directed
spanning tree. The {\it Laplacian matrix} $L=(l_{ij})$ of the
weighted digraph $\mathcal{G}$ is defined by
\begin{equation}
\label{eq202}
l_{ij}=\begin{cases}
-a_{ij}, \qquad i\neq{j},\\
\sum\limits_{k\neq{i}}a_{ik}, \quad i=j.
\end{cases}
\end{equation}
\par
\begin{lemma}
\label{lemma2.1} (Corollary 3 of \cite{scutari2008}) Let
$\mathcal{G}=(\mathcal{V},\mathcal{E},A)$ be a digraph with
Laplacian matrix $L$. If $\mathcal{G}$ is strongly connected, then
$L$ has a simple zero eigenvalue and a positive left-eigenvector
associated to the zero eigenvalue.
\end{lemma}
\begin{definition}
\label{def2.1}
For the high-dimensional Kuramoto model (\ref{eq2}), the {\it
exponential synchronization} is said to be implemented if there
exist $\mu>0$ and a function $\alpha(r(0))$ with respect to the
initial states of the oscillators such that
$$|| r_i(t)-r_j(t)||\leq \alpha(r(0))\  {\rm e}^{-\mu t} \ \ \forall  \ i,\  j=1,2,\cdots,m.$$
\end{definition}
For the high-dimensional Kuramoto model (\ref{eq2}) with identical oscillators ($\Omega_{i}=\Omega$ for each $i=1,2,\cdots,m)$, without loss of generality, we can just consider the dynamics with each $r_{i}$ limited on the unit sphere as follows:
\begin{equation}
\label{eq-6}
\dot{r}_{i}\!=\!k\!\sum\limits_{j=1}^{m}\!a_{ij}(r_{j}\!-\!(r_{i}^{T}\!r_{j})r_{i}),
\ i\!=\!1,2,\cdots,m.
\end{equation}
The first result on the synchronization of $(\ref{eq-6})$ with
digraphs is Theorem 2 of \cite{Lageman-CDC-2016}, which is
rewritten as follows by the terminologies in this paper.
\begin{lemma}
\label{lemma2.2}
 Assume that the digraph ${\mathcal G}$ of the high-dimensional Kuromoto model
 (\ref{eq-6}) has a directed spanning tree and there exists $v\in{\mathbf{R}^{n}}$
 such that $v^\mathrm{T}r_{i}(0)>0$ for every $i=1,2,\cdots,m$.
 Then there exists $\bar{r}\in{\mathbf{R}^n}$ such
 that $\lim\limits_{t\rightarrow +\infty}r_{i}(t)=\bar{r}$ for
 each $i=1,2,\cdots, m$.
\end{lemma}

\begin{remark}
\label{lemma2} In Lemma 2.2, synchronization is implemented on the
hemisphere. But in \cite{Markdahl-Automatica-2018}, a lifting
method is proposed to analyze the global synchronization. In our
recent paper \cite{zhang2018}, the synchronization of (\ref{eq-6})
in the case of digraphs is proved under some stronger conditions
than those in Lemma 2.2 by using a completely different method.

 \end{remark}

\section{Exponential state synchronization}

 In this section, we mainly consider the high-dimensional Kuramoto model (\ref{eq-6}) limited on the unit sphere.\par

Let
\begin{equation}
\label{eq_6}
e_{ij}=1-r_{i}^{T}r_{j}=\frac{1}{2}||r_i-r_j||^2.
\end{equation}
It is easy to see that $e_{ij}=e_{ji}$, $e_{ii}=0$ and
$0\leq{e_{ij}}\leq{2}$ for any $i,j=1,2,\cdots,m$. Obviously,
$e_{ij}=0$ if and only if $r_i=r_j$. So $e_{ij}$ reflects the
state error between $r_i$ and $r_j$. In the following, let us
investigate the  dynamics of all the $e_{ij}$'s. A straightforward
computation shows that
\begin{eqnarray}
\label{eq_08}
\dot{e}_{ij}&=&-r_{j}^{\mathrm T}\dot r_i-r_{i}^{\mathrm T}\dot r_j \nonumber \\
&=&\!-k\!\sum\limits_{l=1}^{m}\!\!a_{il}(\!r_{j}^{\mathrm
T}\!r_{l}\!-\!(\!r_{i}^{\mathrm T}\!r_{\!l})r_{j}^{\mathrm
T}\!r_{\!i})
\!-\!k\!\sum\limits_{l=1}^{m}\!\!a_{\!jl}(r_{i}^{\mathrm
T}\!r_{l}\!-\!(\!r_{j}^{\mathrm T}r_{\!l})r_{i}^{\mathrm
T}\!r_{\!j}).\nonumber \\
\end{eqnarray}
Substituting $r_{i}^{T}r_{j}=1-e_{ij}$ into (\ref{eq_08}) yields
\begin{eqnarray}
\label{eq301}
\dot{e}_{ij} &=& k\sum\limits_{l=1}^{m}\!a_{il}e_{lj}\!-\!k\!\left(\!\sum\limits_{l=1}^{m}a_{il}\!\!\right)\!e_{\!ij} \!-\!k\sum\limits_{l=1}^{m}\!a_{il}e_{li}\nonumber \\
&&
+k\!\left(\!\sum\limits_{l=1}^{m}a_{il}e_{il}\!\!\right)\!e_{\!ij}+ k\sum\limits_{l=1}^{m}a_{jl}e_{li}-k\left(\!\sum\limits_{l=1}^{m}\!a_{jl}\!\!\right)\!e_{ij} \nonumber \\
&&-k\sum\limits_{l=1}^{m}\!\!a_{jl}e_{lj}\!\!+\!k\!\left(\sum\limits_{l=1}^{m}\!\!a_{jl}e_{jl}\!\!\right)\!\!e_{ij}, \ i,j\!=\!\!1,2,\!\cdots\!,m.
\end{eqnarray}
Let $E\!=\!(e_{ij})\!\in{\!\mathbf{R}^{m\times{m}}}$ and
\begin{equation}
\label{eq_8}
\alpha(E\!)\!=\!(\alpha_{1}(E),\alpha_{2}(E),\cdots,\alpha_{m}(E))^{\mathrm T}\!\in{\!\mathbf{R}^{m}},
\end{equation}
where $\alpha_{i}(E)=\sum\limits_{l=1}^{m}a_{il}e_{il}$. Denote by $\Lambda(E)=\mathrm{diag(}\alpha(E))$ the diagonal matrix with the main diagonal elements composed of $\alpha_1(E)$, $\alpha_2(E)$, $\cdots$, $\alpha_m(E)$.
Then we can rewrite (\ref{eq301}) into the compact form described by the Riccati matrix differential equation
\begin{eqnarray}
\label{eq302}
\dot{E}&=&-kLE\!-\!kEL^{\!\mathrm T}\!-k \alpha(E)\mathbf{1}^{\!\mathrm T}\!- \!k\mathbf{1}\alpha^{\!\mathrm T}\!(E)\!
\nonumber \\
&& +k \Lambda (E)E \!+\! kE\Lambda (\!E\!),
\end{eqnarray}
where $L$ is the Laplacian matrix of the digraph.
 \par
In the following, we will use (\ref{eq302}) to investigate whether $E$ converges to zero exponentially. Before our main result on exponential synchronization, we first give a lemma. \par

 \begin{lemma}
\label{lemma4}
Consider a sequence of unit vectors $r_{1},r_{2},\cdots,r_m \in{\mathbf{R}^{n}}$. Let $e_{ij}$ be defined by (\ref{eq_6}) for any $i,j=1,2,\cdots m$.  Then
\begin{equation}
\label{eq303}
 e_{ij}\leq {s(e_{ik_{1}}+e_{k_{1}k_{2}}+\cdots+e_{k_{s-1,j}})}
 \end{equation}
 for any $s\!-\!1$ positive integers $k_1$, $k_2$, $\cdots$,$k_{s-1}$.
\end{lemma}
 \begin{pf}  Since $ \|r_{i}-r_{j}\|^{2}=2e_{ij}$, then  we have
\begin{eqnarray}
\label{eq304}
2e_{ij}&=&\|r_{i}\!-\!r_{j}\|^{2} \nonumber \\
      &\leq& (\|r_{i}\!-\!r_{k_{1}}\!\|\!+\!\|r_{k_{1}}\!-\!r_{k_{2}}\!\|\!+\cdots+\!\|r_{k_{s-1}}\!\!-\!r_{\!j}\|)^{2} \nonumber \\
       &\leq & s(\|r_{i}\!-\!r_{k_{1}}\!\|^{2}\!+\!\|r_{k_{1}}\!\!-\!r_{k_{2}}\!\|^{2}\!+\cdots+\!\|r_{k_{s-1}}\!\!-\!r_{\!j}\|^{2})
       \nonumber \\    &=&
       2s(e_{ik_{1}}+e_{k_{1}k_{2}}+\cdots+e_{k_{s-1},j}),
\end{eqnarray}
where the second inequality comes from the Cauchy inequality.
\hfill $\Box$
 \end{pf}
 \begin{lemma}
\label{lemma3.2}
Assume that the digraph $\mathcal{G}$ has a directed spanning tree with weighted adjacency matrix $A=(a_{ij})\in \mathbf{R}^{m\times m}$. Then \\
{\rm (i)} there exists a constant $c_1>0$ such that
\begin{equation}
\label{eq305}
e_{ij}\leq c_{1}\sum\limits_{p=1}^{m}\sum\limits_{q=1}^ma_{pq}e_{pq},\ \forall \ i,j=1,2,\cdots,m; \end{equation}
{\rm (ii)} there exist constants $\check c>0$ and $\hat c>0$ such that
\begin{equation}
\label{eq306} \hat
c\sum\limits_{i=1}^m\sum\limits_{j=1}^me_{ij}\leq
\sum\limits_{p=1}^{m}\sum\limits_{q=1}^ma_{pq}e_{pq}\leq \check
c\sum\limits_{i=1}^m\sum\limits_{j=1}^me_{ij},
 \end{equation}
that is,
\begin{equation}
\label{eq_14} \hat c \mathbf{1}_m^\mathbf{T}E\mathbf{1}_m\leq
\mathbf{1}_m^\mathbf{T}\alpha(E)\leq \check
c\mathbf{1}_m^\mathbf{T}E\mathbf{1}_m,
 \end{equation}
 where $\alpha(E)$ is defined by (\ref{eq_8}).
\end{lemma}

\begin{pf} (i) Since $\mathcal{G}$ has a directed spanning tree,
there exists a root node $v_{0}$ in $\mathcal{G}$. Then, for any
$i\neq{j}$, there are two directed pathes $v_0k_{1}k_{2}\cdots
k_{s}i$ and $v_0h_{1}h_{2}\cdots h_{t}j$. By Lemma \ref{lemma4},
there exists a constant $c_0>0$
\begin{equation}
e_{i\!j}\!\leq\!
c_0\!(e_{v_0\!k_{1}}\!+e_{k_{1}\!k_{2}}\!+\cdots+\!e_{k_{s}i}\!+\!e_{v_0\!h_{1}}\!+\!e_{h_{1}\!h_{2}}+\cdots+e_{h_{t}\!j}\!)
\nonumber
\end{equation}
for all $i,j=1,2,\cdots,m$. Since all $a_{v_0k_{1}}$,
$a_{k_{1}k_{2}}$, $\cdots$, $a_{k_{s}i}$, $a_{v_0h_{1}}$,
$a_{h_{1}h_{2}}$, $\cdots$, $a_{h_{t}j}$ are positive, there is a
constant $c_{1}>0$ such that
\begin{eqnarray}
\label{eq307}
e_{i\!j}&\!\leq & {c_{1}(a_{v_0k_{1}}e_{v_0k_{1}}+a_{k_{1}k_{2}}e_{k_{1}k_{2}}+\cdots+a_{k_{s}i}e_{k_{s}i})} \nonumber \\
&& +c_{1}(a_{v_0h_{1}}e_{v_0h_{1}}+a_{h_{1}h_{2}}e_{h_{1}h_{2}}+\cdots+a_{h_{t}j}e_{h_{t}j}) \nonumber \\
&\leq&
{c_{1}\sum\limits_{p=1}^{m}\sum\limits_{q=1}^{m}a_{pq}e_{pq}}.
\end{eqnarray}
(ii)  Let $\hat c=\frac{1}{c_{1}m^{2}}$ and $\check
c=\max\limits_{1\leq p,q\leq m}a_{pq}$. Then (\ref{eq306}) follows
from (\ref{eq305}). \hfill $\Box$
\end{pf}
\subsection{The case of strongly connected digraphs}
\begin{theorem}
\label{theorem1} Assume that the digraph ${\mathcal G}$ of the
high-dimensional Kuromoto model  (\ref{eq-6}) is strongly
connected  and there exists $v\in{\mathbf{R}^{n}}$ such that
$v^\mathrm{T}r_{i}(0)>0$ for every $i=1,2,\cdots,m$. Then the
exponential synchronization of  (\ref{eq-6}) is achieved.
\end{theorem}
\begin{pf}  By Lemma \ref{lemma2.2}, we see that the synchronization is achieved
under the conditions of Theorem 3.1. In the rest of the proof, we
use the error dynamics (\ref{eq302}) to prove the exponential
convergence. Since the digraph ${\mathcal G}$ is strongly
connected, by Lemma \ref{lemma2.1}, there is a positive vector
$\beta=(\beta_{1},\cdots,\beta_{m})^{\mathrm
T}\in{\mathbf{R}^{m}}$ satisfying \vspace{-1mm}
\begin{equation}
\label{eqL}
 \beta^{\mathrm T}L=0,\ \ \beta^{T}\mathbf{1}=1.
\end{equation}
 We construct a
total error function as follows:\vspace{-1mm}
$$V(E)=\frac{1}{2}\sum\limits_{i=1}^m\sum\limits_{j=1}^m\beta_i\beta_je_{ij}=\frac{1}{2}\beta^{\mathrm T}E \beta. $$
Let $$\Phi_\eta\!=\!\{E\!=\!(e_{ij})\!\in\mathbf{R}^{m\times m}|\
e_{ij}\!\leq \!\eta,\forall \ i,j\!=\!1,2,\!\cdots\!,m\}$$ and
\vspace{-1mm}
$$\Psi_{\!\eta}\!= \!\{E\!=\!(e_{ij})\!\in\!\mathbf{R}^{m\times m}| \ V\!(E)\!< \!\hat\beta^2\eta/2\},$$
where $0<\eta<1$ and $\hat\beta=\min\limits_{1\leq i\leq
m}\!\!\beta_i>0.$
\\
$~~~~~$ {\bf Claim 1}: $\Psi_\eta\subset \Phi_\eta$.
\par
As a matter of fact, if $E\in \Psi_\eta$, then
$$
 \frac{1}{2}\hat\beta^2e_{ij}\leq \frac{1}{2}\beta_{i}\beta_{j}e_{ij}\leq V(E)<  \frac{1}{2}\hat\beta^2\eta,
$$
which implies that $e_{ij}<\eta$. Thus {\bf Claim 1} is proved.
\par
From (\ref{eqL}), (\ref{eq302}) and {\bf Claim 1}, it follows that
\begin{eqnarray}
\dot{V}(E)&\!=&\!\frac{1}{2}\beta^{\mathrm T}\dot E \beta \nonumber \\
&\!=&\!-k\beta^{\mathrm T} \alpha(E)+k\beta^{T}\Lambda(E) E\beta  \nonumber \\
 &\!=&\!-k\sum\limits_{i=1}^{m}\!\beta_{i}\!\!\sum\limits_{l=1}^{m}\!a_{il}e_{il}
 +k\sum\limits_{i=1}^{m}\sum\limits_{j=1}^{m}\!\beta_{i}
 \!\!\left(\sum\limits_{l=1}^{m}\!a_{il}e_{il}\!\!\right)e_{ij}\beta_j\nonumber \\
 &\!=&\!-k\sum\limits_{i=1}^{m}\sum\limits_{l=1}^{m}(1-\sum\limits_{j=1}^{m}\beta_{j}e_{ij})\beta_{i}a_{il}e_{il}\nonumber\\
\label{eq308} &\leq
&-k(1-\eta)\hat\beta\sum\limits_{i=1}^{m}\sum\limits_{l=1}^{m}a_{il}e_{il},
\ \ \forall \ E\in \Psi_\eta.
\end{eqnarray}
Since the digraph $\mathcal{G}$ is strongly connected, by
(\ref{eq308}) and (\ref{eq306}) of Lemma \ref{lemma3.2}, we have
\begin{eqnarray}
\label{eq309}
\dot{V}(E)&\!\leq &\!-k(1-\eta)\hat\beta \hat c \sum\limits_{i=1}^{m}\sum\limits_{l=1}^{m}e_{il},\nonumber \\
&\!\leq &\!-k(1-\eta)\frac{\hat\beta}{\check\beta^2}\hat c
\sum\limits_{i=1}^{m}
\sum\limits_{l=1}^{m}\beta_i\beta_le_{il},\nonumber \\
&=& -cV(E),\hspace{2cm}\forall \ E\in \Psi_\eta,
\end{eqnarray}
where $\check\beta\!=\!\!\max\limits_{1\leq i\leq m}\!\!\beta_i$
and $c=k(1\!-\!\eta)\frac{\hat\beta}{2\check\beta^2}\hat c\!>\!0.$
By the definition of $\Psi_{\!\!\eta}$ and (\ref{eq309}), we
conclude that $\Psi_{\!\!\eta}$ is a positively invariant set with
respect to (\ref{eq302}), and $V(E)$ converges to zero
exponentially. In details, for any given initial states, there
exists $T$ such that $e_{ij}(T)\in\Psi_{\!\eta }$ since the
synchronization is achieved. By (\ref{eq309}) and the positive
invariance of $\Psi_{\!\!\eta}$, we have that
\begin{equation}
\label{eq_16}
\!V\!(\!E(t))\!\leq \!V\!(\!E(T))\mathrm{e}^{-c(t-\!T)}\!=\!\alpha_{1}(r(0))\mathrm{e}^{-ct},  \ \forall \ t\!>\!T,
\end{equation}
where $\alpha_{1}(r(0))=V(E(T))\mathrm{e}^{cT}$ is
dependent on the initial states.  Moreover, for the time interval $[0,T]$,
there exists $\alpha_2(r(0))$ such that
\begin{equation}
\label{eq_17}
V\!(E(t))\!\leq \!\alpha_2(r(0))\!\leq\!\alpha_2(r(0))\mathrm{e}^{cT}\!\mathrm{e}^{-ct}\!,  \ \forall\ 0\!\leq  \!t\!\leq \!T.
\end{equation}
From (\ref{eq_16}) and (\ref{eq_17}), it follows that there exists $\tilde\alpha(r(0))$ such that
\begin{equation}
\label{eq_18}
V(E(t))\leq \tilde\alpha(r(0))\mathrm{e}^{-ct},  \ \forall \ t \geq 0.
\end{equation}
Further considering
$$
|| r_i(t)-r_j(t))||^2=2e_{ij}(t)\leq \frac{4}{\hat\beta^2}V(E(t)),
$$
we have
$$
|| r_i(t)-r_j(t))||\leq
\frac{2\sqrt{\tilde\alpha(r(0))}}{\hat\beta}\mathrm{e}^{-\frac{c}{2}t}.
$$
Therefore, by Definition \ref{def2.1}, the exponential
synchronization is achieved. \hfill $\Box$
\end{pf}

\begin{remark}
\label{remark2} In the construction of the total error function
$V(E)$, the strong connectedness of ${\mathcal G}$ plays an
important role. From the proof of Theorem \ref{theorem1}, we see
that the crux of the analysis of exponential synchronization is to
demonstrate $\dot V(E)\leq -cV(E)$ for some positive constant $c$
in a neighborhood of the origin $E=0$.
\end{remark}

For the general case of digraphs admitting spanning trees, we will
use all the strongly connected components to construct the total
error function. For the readability of this paper, in the next
subsection we assume that the considered digraph has only two
strongly connected components.
\subsection{The case of digraphs
admitting two strongly connected components. }
\begin{proposition}
\label{lemma5} Suppose that the digraph $\mathcal{G}$ of
(\ref{eq-6}) has a spanning tree and is composed of two strongly
connected components. Then the exponential synchronization is
locally achieved, i.e., there exist a total error function $V(E)$
and a neighborhood $U$ of $E=0$ such that $\dot V(E)\leq -cV(E)$,
$\forall \ E\in U$ for some constant $c>0$.
\end{proposition}
\begin{pf} Let $\mathcal{G}_{1}=(\mathcal{V}_{1},\mathcal{E}_{1})$,
$\mathcal{G}_{2}=(\mathcal{V}_{2},\mathcal{E}_{2})$ be the two
strongly connected components, where $\mathcal{V}_{1}=\{1,2,$
$\cdots,$ $m_1\!\}$,
$\mathcal{V}_{2}\!=\!\{\!{m_1}\!+1,{m_1}\!+2,\cdots,m_1\!+m_2\!\}$
and $m_1\!+m_2\!=\!m$. Without loss of generality, write the
Laplacian matrix as
\begin{equation}
\label{eq_19}
L=\left[\begin{array}{cccc}
L_{1} & 0\\
-A_{21} & L_{2}+D_{2}
\end{array}\right],
\end{equation}
where $A_{21}\ne 0$ is a nonnegative matrix, $D_2$ is a diagnal
nonegative matrix, $L_1\in \mathbf{R}^{m_1\times m_1}$ and $L_2\in
\mathbf{R}^{m_2\times m_2}$ are the Laplacian matrices of
$\mathcal{G}_1$ and  $\mathcal{G}_2$, respectively. Consider the
Lohe model (\ref{eq-6}) and the errors system (\ref{eq302}). Let
$E$ be partitioned as a 2-by-2 block matrix
 \begin{equation}
 \label{eq31 0}
 E=\left[ \begin{array}{cccc}
 E_{11} & E _{12}\\
E^{\mathrm{T}}_{12} & E_{22}
\end{array}\right]
\end{equation}
with $E_{11}\in{\mathbf{R}^{{m_1}\times{{m_1}}}},$
 $E_{22}\in\mathbf{R}^{m_2\times m_2}$.
Since both $\mathcal{G}_1$ and $\mathcal{G}_2$ are strongly
connected, by Lemma \ref{lemma2.1}, there exist positive vectors
$\beta_1=\![\beta_{1\!1},\beta_{12},\!\cdots\!,\beta_{1\!{m_1}}]^{\mathrm{T}}\in
{\mathbf R}^{m_1}$ and
$\beta_2\!\!=\![\beta_{21},\beta_{22},\!\cdots\!,\beta_{2m_2}]^{\mathrm{T}}\in
{\mathbf R}^{m_2}$ such that
\begin{equation}
\label{eq-beta} \beta_1^\mathrm{T}\!L_1\!\!=\!0,\ \
\beta_2^\mathrm{T}\!L_2\!=\!0,\ \
\beta_1^\mathrm{T}\!\mathbf{1}_{m_1}\!\!=\!\beta_2^\mathrm{T}\!\mathbf{1}_{m_2}\!=\!1.
\end{equation}
Let $\beta=[\beta_1^\mathrm{T}\ \
\varepsilon\beta_2^\mathrm{T}]^\mathrm{T}.$ Then
$\beta^\mathrm{T}\mathbf{1}_m=1+\varepsilon$, where
$\varepsilon>0$ is a sufficiently small constant. Construct the
total error function as
$V(E)=\frac{1}{2}\beta^{\mathrm{T}}E\beta$. Then a straightforward
computation shows that
\begin{eqnarray}
\label{eq312}
\dot{V}(E)&\!=&-k\beta^{\mathrm{T}}\!LE\beta-\!k(\!1\!+\!\varepsilon)\beta^{\mathrm{T}}\!\alpha(E)\!+\!k\beta ^{\mathrm{T}}\!\Lambda(\!E) E\beta \nonumber \\
&=&k\varepsilon\beta_{2}^{\mathrm{T}}A_{21}E_{11}\beta_{1}-k\varepsilon\beta_{2}^{\mathrm{T}}D_2E_{12}^{\mathrm{T}}\beta_{1}
\nonumber \\
&&+k\varepsilon^2\beta_{2}^{\mathrm{T}}A_{21}E_{12}\beta_{2}-k\varepsilon^2\beta_{2}^{\mathrm{T}}D_2E_{22}\beta_{2}\nonumber \\
&&-k(\!1\!+\!\varepsilon)(\beta_{\!1}^{\!\mathrm{T}}\!\alpha_{\!1}\!(\!E)\!+\!\varepsilon\beta_{2}^{\!\mathrm{T}}\!\alpha_{2}(\!E)\!)\!+\!k\beta ^{\mathrm{T}}\!\Lambda(\!E) \!E\beta\nonumber \\
&\leq&k\varepsilon\beta_{2}^{\mathrm{T}}A_{21}E_{11}\beta_{1}+k\varepsilon^2\beta_{2}^{\mathrm{T}}A_{21}E_{12}\beta_{2}\nonumber\\ &&-k\beta_{\!1}^{\!\mathrm{T}}\!\alpha_{\!1}\!(\!E)\!-\!k\varepsilon(\beta_{\!1}^{\!\mathrm{T}}\!\alpha_{\!1}\!(\!E)\!+\!\beta_{2}^{\!\mathrm{T}}\!\alpha_{2}(\!E))\!+\!\gamma(\!E),
\end{eqnarray}
where $\gamma(E)=\!k\beta ^{\mathrm{T}}\!\Lambda(\!E) \!E\beta$
contains all the higher-order terms with respect to $E$. Letting
$$\check\beta_i\!=\!\!\max\limits_{1\leq j\leq m_i}\!\!\beta_{ij},\ \
\hat\beta_i\!=\!\!\min\limits_{1\leq j\leq m_i}\!\!\beta_{ij},\ \
\check a \!=\!\max(A_{21}),$$ by (\ref{eq312}) and
(\ref{eq-beta}), we have
\begin{eqnarray}
\label{eq_25} \dot{V}(E)&\leq& (k\varepsilon\check a \check
\beta_1)\mathbf{1}_{m_1}^{\mathrm{T}}E_{11}\mathbf{1}_{m_1}+(k\varepsilon^2\check\beta_{2}\check
a)\mathbf{1}_{m_1}^{\mathrm{T}}E_{12}\mathbf{1}_{m_2}\nonumber\\
&&-(k\hat\beta_{\!1})
\mathbf{1}_{m_1}^{\mathrm{T}}\alpha_{\!1}\!(\!E)\!-\!(k\varepsilon\hat\beta_{\!1}\hat\beta_{2})
\mathbf{1}_m^{\mathrm{T}}\alpha(\!E)\!+\!\gamma(\!E).
\end{eqnarray}
Since ${\mathcal G}_1$ is strongly connected and ${\mathcal G}$
has a spanning tree, by (\ref{eq306}) in Lemma \ref{lemma3.2} and
(\ref{eq_25}), we have
\begin{eqnarray}
\label{eq_26}
\dot{V}(E)&\leq& k(\varepsilon\check\beta_{2}\check a \check \beta_1- \hat\beta_{\!1}\hat c_{1})
\mathbf{1}_{m_1}^{\mathrm{T}}E_{11}\mathbf{1}_{m_1}\nonumber\\
& & +k\varepsilon(\varepsilon\check\beta_{2}^2\check
a-\hat\beta_{\!1}\hat\beta_{2}\hat
c)\mathbf{1}_m^{\mathrm{T}}E_{}\mathbf{1}_{m}\!+\!\gamma(\!E),
\end{eqnarray}
where $\hat c_1$ and $\hat c$ are determined by Lemma
\ref{lemma3.2}. By (\ref{eq_26}), as $\varepsilon$ is sufficiently
small, there is a constant $\tilde c>0$ such that
\begin{eqnarray}
\label{eq_27}
\dot{V}(E)&\leq& -\tilde c \ \mathbf{1}_m^{\mathrm{T}}E_{}\mathbf{1}_{m}\!+\!\gamma(\!E)
\end{eqnarray}
Since $\gamma(E)$ is composed of all the higher-order terms, there exists a neighborhood $U$ of $E=0$ such that
\begin{equation}
\label{eq_28} \dot{V}\!(E)\!\leq\! -\frac{\tilde c}{2}
\mathbf{1}_m^{\mathrm{T}}E\mathbf{1}_{m}\!\leq \!-\frac{\tilde
c}{2\hat \beta^2}\beta^{\mathrm{T}}\!E\beta\!=\!-cV\!(E),  \ \
\forall\ E\in U, \nonumber
\end{equation}
where $c=\tilde c/\hat \beta^2>0$. \hfill $\Box$
\end{pf}
\subsection{The case of digraphs admitting spanning trees. }
\begin{theorem}
\label{theorem2} If the digraph $\mathcal{G}$ of the
high-dimensional Kuramoto model (\ref{eq-6}) has a spanning tree,
then the exponential synchronization is locally achieved.
\end{theorem}
%
\begin{pf} Denote all the strongly connected components of $\mathcal{G}$ by
$\mathcal{G}_1$, $\mathcal{G}_2$, $\cdots$, $\mathcal{G}_\mu$, respectively.
Let the number of the nodes of $\mathcal{G}_i$ be $m_i$ for each $i=1,2,\cdots,\mu$.
Without loss of generality, assume that the Laplacian matrix $L$ has the form as follows:
\begin{equation}
\label{eq20}
L
=\!\left[\!\!\begin{array}{cccc}
L_1& &\ & \\
-\!A_{21} & L_2\!+\!\!D_2 & \\
\vdots  &\!\!\!\!\!\!\ddots  &\ddots \\
-\!A_{\mu 1} & \!\!\!\!\!\!\!\cdots & \!\!\!\!\!\!\!\!\!-\!A_{\mu,\mu-1 } &L_\mu\!\!+\!\!D_\mu
\end{array}\!\!\right]\!\!,
\end{equation}
where $L_i$ is the Laplacian matrix of $\mathcal{G}_i$, each
$A_{ij}\in \mathbf{R}^{m_i\!\times \!m_j}$ is a nonnegative
matrix, and every $D_{i}\in \mathbf{R}^{m_i\!\times \!m_i}$ is a
nonnegative diagonal matrix.  
Corresponding to (\ref{eq20}),  we partition $E$ and $\alpha(E)$
as
\begin{equation}
\label{eq21} E\!\!=\!\! \left[\!\!\begin{array}{ccc}
E_{11}& \cdots &E_{1\mu}\\
\vdots & \ddots &  \vdots\\
E_{\mu1}& \cdots & E_{\mu\mu}\\
\end{array}\!\!\right],\ \
\alpha(E)=\left[\!\!\begin{array}{c}
\alpha_1(E)\\
\vdots \\
\alpha_\mu(E) \\
\end{array}\!\!\right],
\end{equation}
respectively, where
$E_{ij}=E_{ij}^\mathrm{T}\in{\mathbf{R}^{m_i\times m_j}}$ for each
$i,j=1,2\cdots,\mu$.
Since $\mathcal{G}_i$ is strongly connected, by Lemma
\ref{lemma2.1}, there exists a positive vector $\beta_i\in
{\mathbf R}^{m_i}$ such that
\begin{equation}
\label{eq_23} \beta_i^{\mathrm{T}}L_i=0, \ \
\beta_i^{\mathrm{T}}\mathbf{1}_{m_i}=1
\end{equation}
for each $i=1,2,\cdots,\mu$. Construct the total error function
$V\!(E)\!=\!\frac{1}{2}\beta^{\!\mathrm{T}}\!E\beta$ with
\begin{equation}
\label{eq_24} \beta =[\beta_1^{\mathrm{T}},\ \
\varepsilon\beta_2^{\mathrm{T}},\ \
\varepsilon^2\beta_3^{\mathrm{T}},\cdots,
\varepsilon^{\mu-\!1}\beta_\mu^{\mathrm{T}}]^{\mathrm{T}}.
\end{equation} Then
$\beta^\mathrm{T}\mathbf{1}_m=1+\varepsilon+\varepsilon^2+\cdots+\varepsilon^{\mu-1}$.
A straightforward computation shows that
\begin{eqnarray}
\label{eq_34}
\dot{V}(E)&\!=&-k\beta^{\mathrm{T}}\!LE\beta-\!k\beta^{\mathrm{T}}\!\alpha(E)\mathbf{1}_m^\mathrm{T}\beta+\!k\beta ^{\mathrm{T}}\!\Lambda(\!E) E\beta \nonumber \\
&\leq &k\sum_{i=2}^\mu\sum_{j=1}^{i-1}\sum_{l=1}^\mu\varepsilon^{i+l-2}\beta_{i}^{\mathrm{T}}A_{ij}E_{jl}\beta_{l}\nonumber \\
&& -k\sum_{i=1}^\mu\varepsilon^{i-1}\sum_{j=1}^\mu\varepsilon^{j-1} \beta_{j}^{\mathrm{T}}\alpha_j(E)+\!\gamma(\!E),
\end{eqnarray}
where $\gamma(E)=\!k\beta ^{\mathrm{T}}\!\Lambda(\!E) \!E\beta$
contains all the higher-order terms with respect to $E$. Denote by
$\mathcal{\tilde G}_i$ the subgraph composed of the first $i$
strongly connected components of $\mathcal{G}$. Let $\tilde
m_i=m_1+m_2+\cdots+m_i$ and
\begin{equation}
\tilde E_i\!\!=\!\! \left[\!\!\begin{array}{ccc}
E_{11}& \cdots &E_{1i}\\
\vdots & \ddots &  \vdots\\
E_{i1}& \cdots & E_{ii}\\
\end{array}\!\!\right],\ \
\tilde\alpha_i(E)=\left[\!\!\begin{array}{c}
\alpha_1(E)\\
\vdots \\
\alpha_i(E) \\
\end{array}\!\!\right].
\end{equation}
By the non-negativity of  $A_{ij}$ and  $\beta_i$, there exists a
constant $d_1>0$ such that
\begin{eqnarray}
\label{eq_35} &\beta_{i}^{\mathrm{T}}A_{ij}E_{jl}\beta_{l}\leq
&d_{1}\mathbf{1}_{m_j}^{\mathrm{T}}E_{jl}\mathbf{1}_{m_i}\leq
d_{1}\mathbf{1}_{\tilde m_{j\vee l}}^{\mathrm{T}}\tilde E_{j\vee
l}\mathbf{1}_{\tilde m_{j\vee l}},
\end{eqnarray}
where ${j\vee l}=\max\{j,l\}$. Moreover,  by (\ref{eq_14}) in
Lemma \ref{lemma3.2}  and the positivity of $\beta_i$'s, it is
easily seen that there exist positive constants $d_2$ and $d_3$
both independent of $\varepsilon$ such that
\begin{eqnarray}
\label{eq_36}
&& (1\!+\!\varepsilon \!+\!\cdots\!+\!\varepsilon^{\mu-\!1}\!)(\beta_{\!1}\!^{\mathrm{T}}\alpha_{\!1}\!+\!\varepsilon\beta_2^{\mathrm{T}}\alpha_2\!+\!\cdots\!+\!\varepsilon^{\mu-\!1}\!\beta_{\!\mu}\!^{\mathrm{T}}\alpha_{\!\mu}\!)\nonumber\\
&\geq &\beta_{\!1}\!^{\!\mathrm{T}}\!\alpha_{\!1}\!+
\!\varepsilon(\beta_{\!1}\!^{\!\mathrm{T}}
\!\alpha_{\!1}\!+\!\beta_{\!2}\!^{\!\mathrm{T}}\!\alpha_{\!2}\!)\!+\!\cdots\!+
\!\varepsilon^{\mu-\!1}\!(\beta_{\!1}\!^{\!\mathrm{T}}
\!\alpha_{\!1}\!+\!\cdots\!+\!\beta_{\!\mu}^{\!\mathrm{T}}\!\alpha_{\!\mu}\!)\nonumber\\
&\geq &d_{2}\sum_{s=1}^\mu\varepsilon^{s-1}\mathbf{1}_{\tilde
m_s}^{\mathrm{T}}\!\tilde\alpha_{s} \!\nonumber\\
&\geq &d_{3}\sum_{s=1}^\mu\varepsilon^{s-1}\mathbf{1}_{\tilde
m_{s}}^{\mathrm{T}} \tilde E_{s}\mathbf{1}_{\tilde m_{s}}.
 \end{eqnarray}
Applying (\ref{eq_35}) and (\ref{eq_36}) to (\ref{eq_34}) yields
\begin{eqnarray}
\label{eq_37} \dot{V}(E) &\leq
&kd_{1}\sum_{i=2}^\mu\sum_{j=1}^{i-1}\sum_{l=1}^\mu\varepsilon^{i+l-2}
\mathbf{1}_{\tilde m_{j\vee l}}^{\mathrm{T}}\tilde E_{j\vee
l}\mathbf{1}_{\tilde m_{j\vee l}}\nonumber \\ &&
-kd_{3}\sum_{s=1}^\mu \varepsilon^{s-1}\mathbf{1}_{\tilde
m_s}^{\mathrm{T}}\tilde E_s\mathbf{1}_{\tilde m_s}+\!\gamma(\!E).
\end{eqnarray}
From (\ref{eq_37}), we see that
\begin{equation}
\label{eq_38}
i\!+\!l\!-\!2\!\geq \!i\!-\!1\!\geq \!j,\ \ i\!+\!l\!-\!2\!\geq \!l.
\end{equation}
It follows from (\ref{eq_38}) that $i\!+\!l\!-\!2\geq j\vee l$. So, by (\ref{eq_37}),we have
\begin{eqnarray}
\label{eq_39}
\dot{V}\!(\!E)
&\!\leq & \!-k\frac{d_{3}}{2}\sum_{s=1}^\mu \!\varepsilon^{s-\!1}\!\mathbf{1}_{\tilde m_s}^{\!\mathrm{T}}\!\!\tilde E_s\mathbf{1}_{\tilde m_s}\!\!\!+\!\gamma(\!E)\nonumber \\ &&\leq \!\!-k\frac{d_{3}}{2}\varepsilon^{\mu-\!1}\mathbf{1}_{m}^{\mathrm{T}}E\mathbf{1}_{m}\!+\!\gamma(\!E)
\end{eqnarray}
for sufficiently small $\varepsilon>0$. So there exist a constant
$c>0$ and a neighborhood $U$ of the origin $E=0$ such that $\dot
V(E)\leq -cV(E)$, which implies that $V(E)$ converges to zero
exponentially.\hfill $\Box$
\end{pf}
Combining Lemma \ref{lemma2.2} and Theorem \ref{theorem2} yields
the following result:
\begin{theorem}
\label{theorem3} Assume that the digraph ${\mathcal G}$ of the
high-dimensional Kuromoto model  (\ref{eq-6}) has a directed
spanning tree and there exists $v\in{\mathbf{R}^{n}}$ such that
$v^\mathrm{T}r_{i}(0)>0$ for every $i=1,2,\cdots,m$. Then the
exponential synchronization of (\ref{eq-6}) is achieved.
\end{theorem}
Assume that the high-dimensional Kuramoto model with identical
oscillators has the general form as follows
\begin{equation}
\label{eqOmega} \dot{r}_{i}\!=\Omega
r_i+\!k\!\sum\limits_{j=1}^{m}\!a_{ij}(r_{j}\!-\!(r_{i}^{T}\!r_{j})r_{i}),
\ i\!=\!1,2,\cdots,m,
\end{equation}
where $\Omega$ is a skew-symmetric matrix. Then, by
\cite{zhu2013}, a system transformation can be adopted to
transform system (\ref{eqOmega}) to the form of (\ref{eq-6}).
Therefore, the conclusion of Theorem \ref{theorem3} still holds
for system (\ref{eqOmega}).
 \section{Simulations}
 In this section, we give some simulations to show the exponential synchronization.
Consider the high-dimensional Kuramoto model (\ref{eq-6}) with
$n=3$ and $m=12$. The digraph is shown in Fig. 1, which has 3
strongly connected components.
\begin{figure}[t]
      \centering
      \includegraphics[height=3.6cm,width=6.2cm]{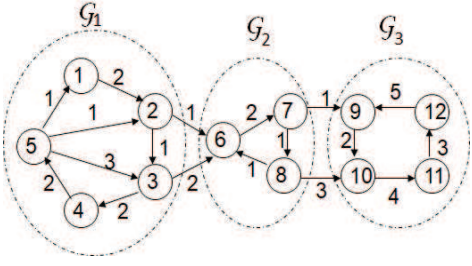}
          \caption{The digraph of the high-dimensional Kuramoto model.}
      \label{figurelabel}
   \end{figure}
From Fig. 1, we obtain the adjacency matrix as follows: {\small
\setlength{\arraycolsep}{2.4pt}
\begin{equation}
\label{eq_40}
  A=(a_{ij})= \left[\begin{array}{cccccccccccc}
0 & 0 & 0 & 0 & 1& 0 & 0 & 0 & 0 &0 & 0 & 0 \vspace{-1.5mm}\\
2 & 0 & 0 & 0 & 1& 0 & 0 & 0 & 0 &0 & 0 & 0 \vspace{-1.5mm}\\
0 & 1 & 0 & 0 & 3& 0 & 0 & 0 & 0 &0 & 0 & 0 \vspace{-1.5mm}\\
0 & 0  & 2& 0 & 0& 0 & 0 & 0 & 0 &0 & 0 & 0 \vspace{-1.5mm}\\
0 & 0 & 0 & 2 & 0& 0 & 0 & 0 & 0 &0 & 0 & 0 \vspace{-1.5mm}\\
0 & 1 & 2 & 0 & 0& 0 & 0 & 1 & 0 &0 & 0 & 0 \vspace{-1.5mm}\\
0 & 0  & 0  & 0 & 0& 2& 0 & 0 & 0 &0 & 0 & 0 \vspace{-1.5mm}\\
0 & 0  & 0  & 0 & 0& 0 & 1 &0 & 0 &0 & 0 & 0 \vspace{-1.5mm}\\
0 & 0  & 0  & 0 & 0& 0 & 1 & 0& 0 &0 & 0 & 5 \vspace{-1.5mm}\\
0 & 0  & 0  & 0 & 0& 0 & 0 & 3& 2 &0 & 0 & 0 \vspace{-1.5mm}\\
0 & 0  & 0  & 0 & 0& 0 & 0 & 0 & 0  &4& 0 & 0 \vspace{-1.5mm}\\
0 & 0  & 0  & 0 & 0& 0 & 0 & 0 & 0  &0 & 3& 0
\end{array}\right].
\end{equation}
}
Then, the dynamics with each $r_{i}$ limited on the unit sphere is
\begin{equation}
\dot{r}_{i}\!=\!k\!\sum\limits_{j=1}^{12}\!a_{ij}(r_{j}\!-\!(r_{i}^{T}\!r_{j})r_{i}),
\ i\!=\!1,2,\cdots,12.
\end{equation}
From Fig. 1, the Laplacian matrices of $\mathcal{G}_1$, $\mathcal{G}_2$ and $\mathcal{G}_3$ are
$$
{ \setlength{\arraycolsep}{0.8pt}
L_{\!1}\!\!=\!\!\!\left[\!\!\begin{array}{ccccc}
1 & 0 & 0 & 0 & \!-1 \vspace{-1mm}\\
-2 & 3 & 0 & 0 & \!-1 \vspace{-1mm}\\
0 & \!-1 & 4 & 0 & \!-3 \vspace{-1mm}\\
0 & 0  & \!-2& 2 & 0 \vspace{-1mm}\\
0 & 0 & 0 & \!-2 & 2
\end{array}\!\!\right]\!\!,
L_{\!2}\!\!=\!\!\!\left[\!\!\begin{array}{ccc}
 1 & 0 & -1  \vspace{-1mm}\\
 -2 & 2 & 0  \vspace{-1mm}\\
 0 & -1 & 1
\end{array}\!\!\right]\!\!,
L_{\!3}\!\!=\!\!\!\left[\!\!\begin{array}{cccc}
 5 &0 & 0 & \!-5 \vspace{-1mm}\\
 -2 &2 & 0 & 0 \vspace{-1mm}\\
 0  &\!-4& 4 & 0 \vspace{-1mm}\\
 0  &0& \!-3 & 3
\end{array}\!\!\right]\!\!,
}$$
respectively. A straightforward computation shows that
$$
\beta_1^\mathrm{T}=[0.1111,\  0.0556,\  0.1667,\ 0.3333,\  0.3333],
\vspace{-3mm}
$$
$$
\beta_{\!2}^{\!\mathrm{T}}\!\!\!=\!\![0.4, \   0.2,\ 0.4],\
 \beta_{\!3}^{\!\mathrm{T}}\!\!\!=\!\![0.1558,\ 0.3896,\ 0.1948, 0.2597].$$
 Construct the total error function $V\!(E)\!=\!\frac{1}{2}\beta^{\!\mathrm{T}}\!E\beta$ with
\begin{equation}
\label{eq_24} \beta =[\beta_1^{\mathrm{T}},\ \
\varepsilon\beta_2^{\mathrm{T}},\ \
\varepsilon^2\beta_3^{\mathrm{T}}]^{\mathrm{T}}. \end{equation}
Fig.2 shows that the time response curves of agent's are
synchronized. Fig.3 shows that all the trajectories of  the
oscillators converge to the same point on the unit sphere.

Consider the high-dimensional Kuramoto model with the form as
follows:
\begin{equation}
\dot{r}_{i}\!=\Omega
r_i+\!k\!\sum\limits_{j=1}^{12}\!a_{ij}(r_{j}\!-\!(r_{i}^{T}\!r_{j})r_{i}),
\ i\!=\!1,2,\cdots,1\!2,
\end{equation}
where
$$
\Omega=\left[\begin{array}{ccc}
                     0 &1& -2\\
                     -1& 0& -1\\
                     2 & 1& 0\end{array}\right].
$$

Fig.4 shows that the dynamical synchronization is achieved. Fig.5
shows that the trajectories converge to a periodic orbit. Finally,
we display the total error function $V$. In Fig.6, it is shown
that the total error function tends to zero exponentially as
$t\rightarrow\infty$.
\begin{figure}[t]
      \centering
      \includegraphics[height=6cm,width=7cm]{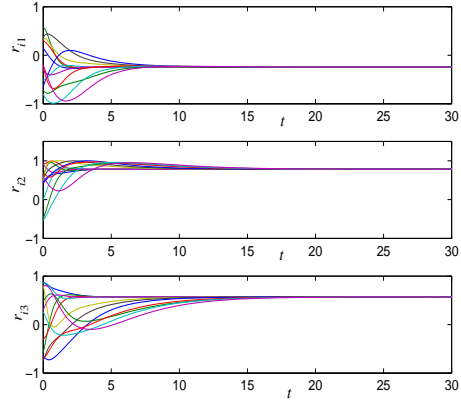}
          \caption{The time response curves as $k_i=1$.}
      \label{figurelabel}
   \end{figure}

\begin{figure}[t]
      \centering
      \includegraphics[height=6cm,width=7cm]{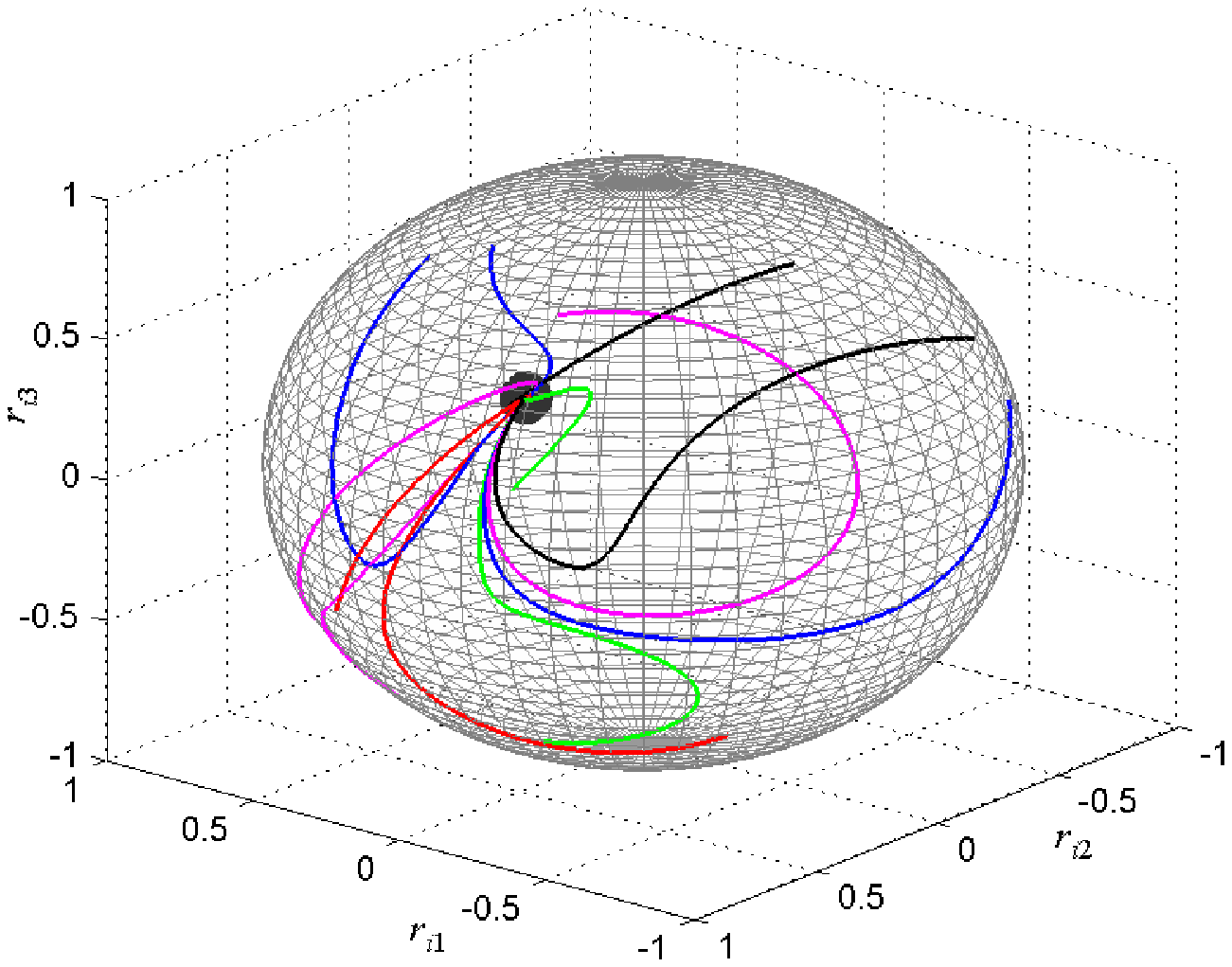}
          \caption{The trajectories of the oscillators as $k_i=1$.}
      \label{figurelabel}
   \end{figure}

\begin{figure}[t]
      \centering
      \includegraphics[height=6cm,width=7cm]{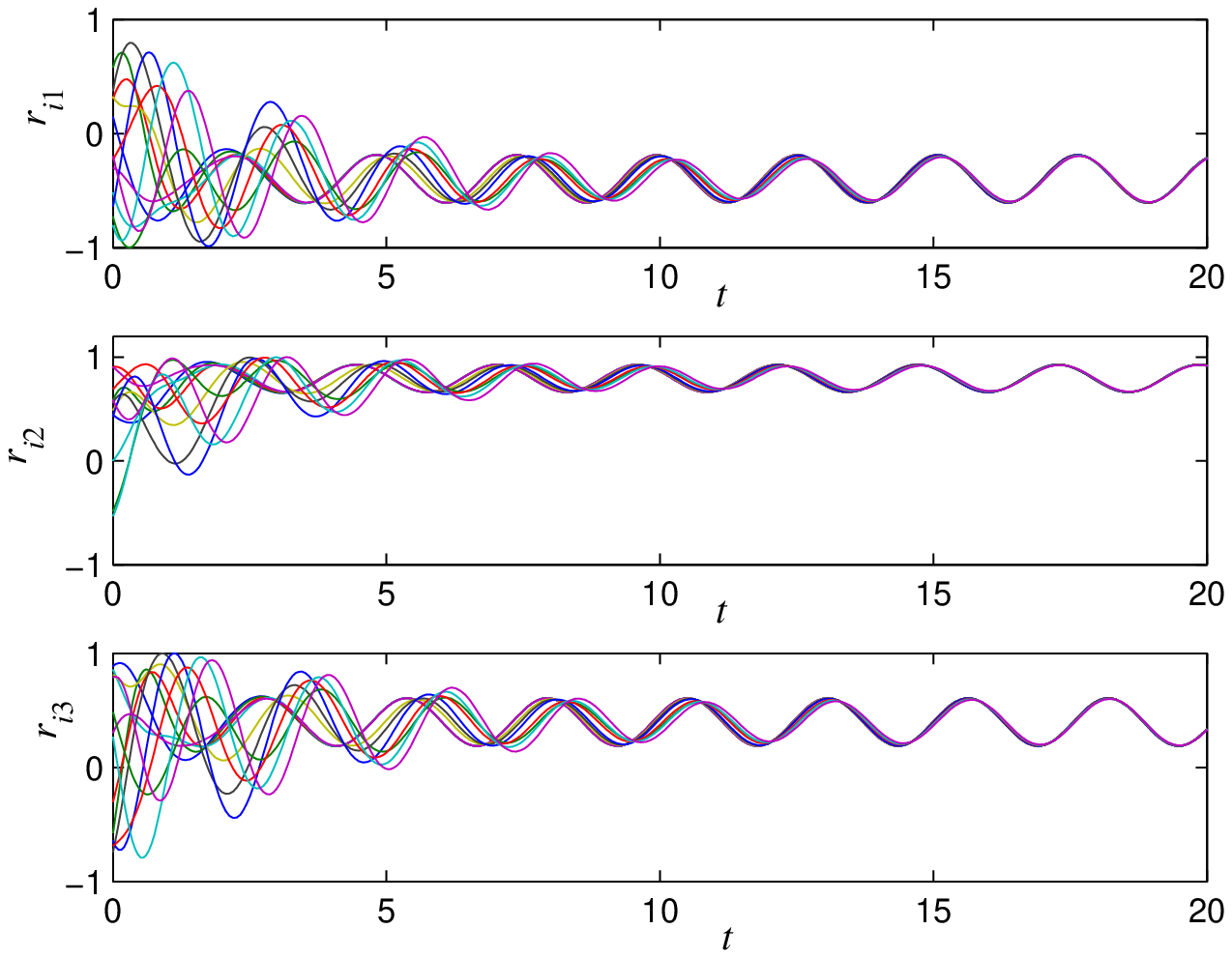}
          \caption{The time response curves as $k_i=1$.}
      \label{figurelabel}
   \end{figure}

\begin{figure}[t]
      \centering
      \includegraphics[height=6cm,width=7cm]{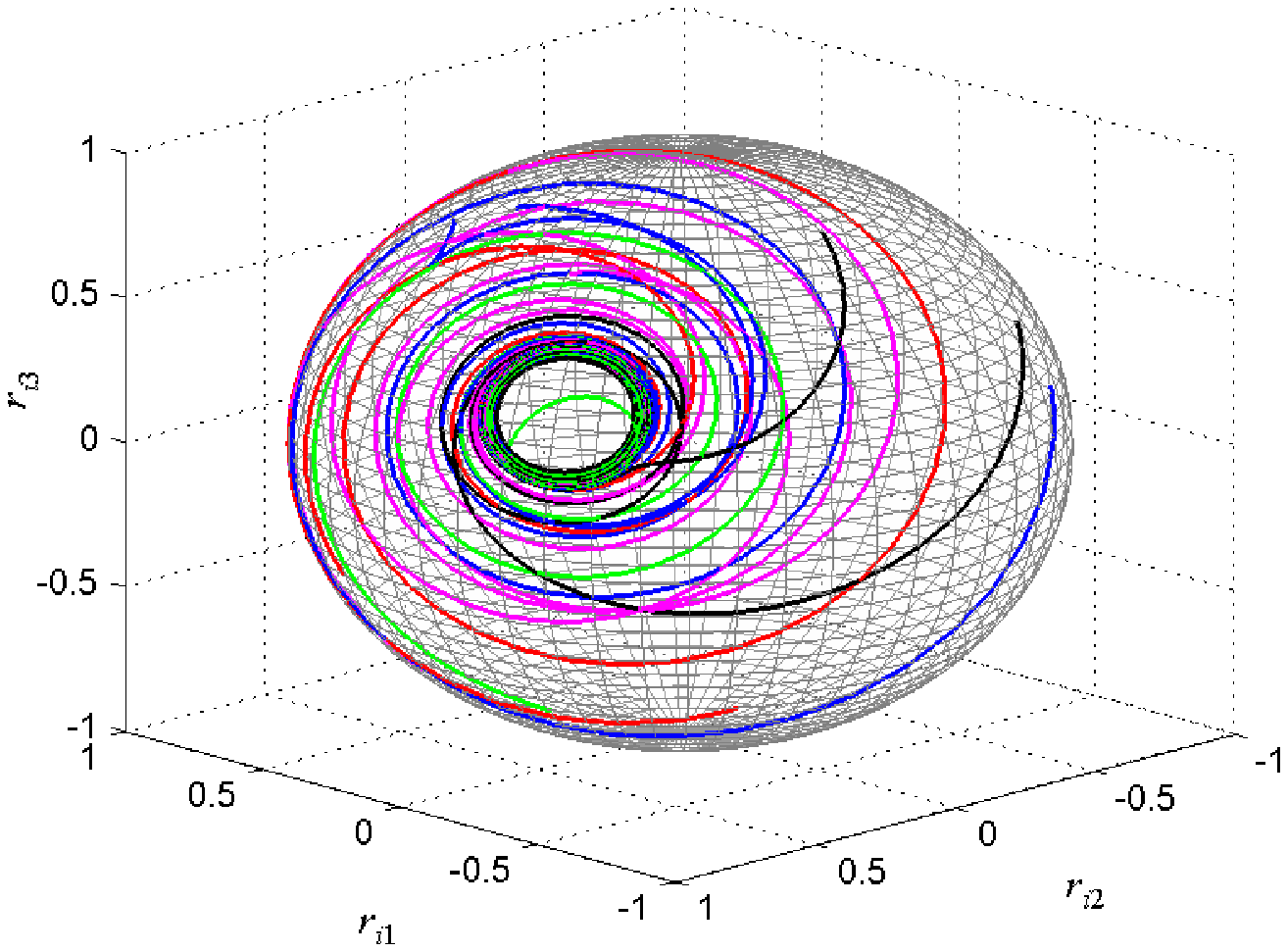}
          \caption{The trajectories of the oscillators as $k_i=1$.}
      \label{figurelabel}
   \end{figure}
\begin{figure}[t]
      \centering
      \includegraphics[height=6cm,width=7cm]{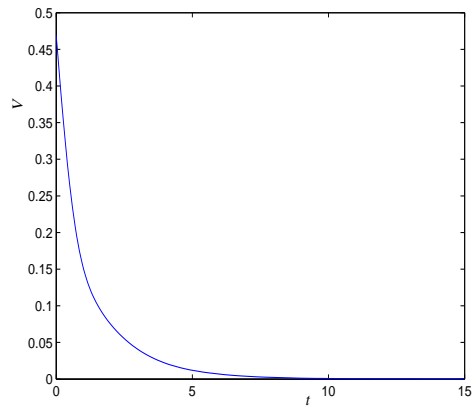}
          \caption{The time response curves of the total error function.}
      \label{figurelabel}
   \end{figure}
\section{Conclusions}
The exponential synchronization has been proved for the
high-dimensional Kuramoto model with identical oscillators under
the digraphs admitting spanning trees. The error dynamics is
described by a matrix Riccati differential equation and a total
error function is constructed for the analysis of the exponential
synchronization. In our future work, the exponential
synchronization will be investigate for some generalized
high-dimensional Kuramoto model. For the high-dimensional Kuramoto
model with non-identical oscillators, practical synchronization
will be addressed. How to generalize the frequency synchronization
into the high-dimensional spaces is also an interesting issue.
\bibliographystyle{abbrv}

\end{document}